\documentclass[11pt]{article}
\usepackage{amssymb}
\newtheorem{theorem}{Theorem}

\newtheorem{proposition}{Proposition}

\newtheorem{remark}{Remark}

\def\defi{\stackrel{{\scriptscriptstyle \Delta}}{=}}

\def\a{\alpha}

\def\o{\omega}

\def\F{{\cal F}}
\def\w{\widehat}

\def\Arg{{\rm Arg\,}}

\def\Re{{\rm Re\,}}

\def\R{{\bf R}}

\def\L{L}

\def\b{\beta}

\def\C{{\bf C}}
\def\W{{\cal W}^*}
\def\W{{\cal W}}
\def\ww{\widetilde}

\def\D{{\Delta}}
\def\p{\partial}

\def\L{{\cal L}}

\newcommand{\be}{\begin{equation}}
\newcommand{\ee}{\end{equation}}
\newcommand{\bd}{\begin{displaymath}}
\newcommand{\ed}{\end{displaymath}}
\newcommand{\ba}{\begin{array}{ll}}
\newcommand{\ea}{\end{array}}
\newcommand{\baa}{\begin{eqnarray}}
\newcommand{\eaa}{\end{eqnarray}}
\newcommand{\baaa}{\begin{eqnarray*}}
\newcommand{\eaaa}{\end{eqnarray*}}
\font\sm=cmr10

\def\a{\alpha}
\def\D{\cal C}
\title{
Parabolic equations with the second order Cauchy conditions on
the boundary}
\author{
Nikolai Dokuchaev\\ {\sm  Department of Mathematics, Trent
University, Ontario, Canada}}
 
\begin{document}
 \vspace{-0.5cm}      \maketitle
\begin{abstract}
The paper studies some ill-posed boundary value problems on
semi-plane for parabolic equations with homogenuous Cauchy
condition at initial time and with the second order Cauchy
condition on the boundary of the semi-plane. A class of inputs
that allows some regularity is suggested and described explicitly
in frequency domain. This class is everywhere dense in the space
of square integrable functions.
\end{abstract}
Parabolic equations such as heat equations have fundamental
significance for natural sciences, and various boundary value
problems for them were widely studied  including well-posed
problems as well as the so-called  ill-posed problems that are
often significant for applications. The present paper introduces
and investigates a special boundary value problem on semi-plane
for parabolic equations with homogenuous Cauchy condition at
initial time and with second order Cauchy condition on the
boundary of the semi-plane. The problem is ill-posed. A set of
solvability, or a class of inputs that allows some regularity in a
form of prior energy type estimates is suggested and described
explicitly in frequency domain. This class is everywhere dense in
the class of $L_2$-integrable functions. This result looks
counterintuitive, since these boundary conditions are unusual;
solvability of this boundary value problem for a wider class of
inputs  is inconsistent with basic theory.
\section{The problem setting}
Let us consider the following boundary value problem \baa
&&a\frac{\p u}{\p t}(x,t)=\frac{\p^2 u}{\p x^2}(x,t)+b\frac{\p
u}{\p x}(x,t)+c u(x,t)+f(x,t),\nonumber\\
&&\hphantom{a}u(x,0)\equiv 0, \nonumber\\
&&\hphantom{a}u(0,t)\equiv g_0(t),\quad \frac{\p u}{\p
x}(0,t)\equiv g_1(t). \label{parab} \eaa
 Here $x>0$, $t>0$, and $a>0,b,c\in\R$ are constants, $g_k\in L_2(0,+\infty)$, $k=1,2$, and $f$ is a
 measurable function such that
 $\int_0^ydx\int_0^{\infty}|f(x,t)|^2dt<+\infty$ for all $y>0$.
 \par
 This problem is ill-posed (see Tikhonov and Arsenin (1977)).
 \par
 Let $\mu\defi b^2/4-c$. We assume that $\mu>0$. Note that this assumtion does not reduce generality
 for the cases when we are interested in solution on a finite time interval, since we can rewrite
 the parabolic equation as the one with $c$ replaced by $c-M$ for any $M>0$ and $g_k(t)$ replaced by
 $e^{-Mt}g_k(t)$; the solution $u_M$ of the new equation
 related to the solution $u$ of the old one as
 $u_M(x,t)=e^{-Mt}u(x,t)$.
\subsection*{Definitions and special functions}
Let $\R^+\defi[0,+\infty)$, $\C^+\defi\{z\in\C:\ \Re z>  0\}$. For
$v\in L_2(\R)$, we denote by $\F v$ and $\L v$ the Fourier and the
Laplace  transforms respectively
 \baa\label{U}
V(i\o)=(\F v)(i\o)\defi\frac{1}{\sqrt{2\pi}}\int_{\R}e^{-i\o
t}v(t)dt, \quad\o\in\R,\eaa
 \baa\label{Up}
V(p)=(\L v)(p)\defi\frac{1}{\sqrt{2\pi}}\int_{0}^{\infty}e^{-p
t}v(t)dt, \quad p\in\C^+. \eaa
\par
Let  $H^r$ be the Hardy space of holomorphic on $\C^+$ functions
$h(p)$ with finite norm
$\|h\|_{H^r}=\sup_{k>0}\|h(k+i\o)\|_{L_r(\R)}$, $r\in[1,+\infty]$
(see, e.g., Duren (1970)).
\par
For $y>0$, let  $\W(y)$ be the Banach space of the functions
$u:(0,y)\times\R^+\to\R$ with the finite norm \baaa \|u\|_{\W(y)}
\defi
\sup_{x\in(0,y)}\biggl(\|u(x,\cdot)\|_{L_2(\R^+)}+\Bigl\|\frac{\p
u}{\p x} (x,\cdot)\Bigr\|_{L_2(\R^+)}&+&\Bigl\|\frac{\p^2 u}{\p
x^2} (x,\cdot)\Bigr\|_{L_2(\R^+)}\\&+&\Bigl\|\frac{\p u}{\p t}
(x,\cdot)\Bigr\|_{L_2(\R^+)}\biggr). \eaaa
The class $\W(y)$ is such that all the equations presented in problem
(\ref{parab}) are well defined for any $u\in\W(y)$ and in the domain $(0,y)\times\R^+$.
For instance,
If $v\in \W(y)$, then, for any $t_*>0$, we have that
$v|_{[0,y]\times[0,t_*]}\in C([0,t_*],L_2(0,y))$ as a function of
$t\in[0,t_*]$. Hence the initial condition at time $t=0$ is well
defined as an equality in $L_2([0,y])$.  Further,
we have that  $v|_{[0,y]\times\R^+}\in C([0,y],L_2(\R^+))$ and
$\frac{\p v}{\p x}\Bigr|_{[0,y]\times\R^+}\in
C([0,y],L_2(\R^+))$ as functions of $x\in[0,y]$. Hence the
functions $v(0,t)$, $
\frac{dv}{dx}(x,t)|_{x=0}$ are well defined as elements of
$L_2(\R^+)$, and the boundary value conditions at $x=0$ are well
defined as equalities in $L_2(\R^+)$.
\subsection*{Special smoothing kernel}
 Let us introduce the set of the following special function:
\baa K(p)=K_{\a,\b,q}(p)\defi e^{-\a (p+\b)^q},\quad p\in\C^+.
\label{K}\eaa Here $\a>0$, $\b>0$ are reals, and
$q\in(\frac{1}{2},1)$ is a rational number.
 We mean the branch of $(p+\b)^q$
such that its argument is $q\Arg(p+\b)$, where $\Arg z\in
(-\pi,\pi]$ denotes the principal value of the argument of
$z\in\C$.
\par
 The functions $K_{\a,\b,q}(p)$ are holomorphic in $\C^+$, and $$ \ln|K(p)|=-\Re(\a
(p+\b)^q)=-\a|p+\b|^q\cos [q\Arg (p+\b)].$$  In addition, there
exists $M=M(\b,q)>0$ such that $\cos[q\Arg(p+\b)]>M$ for all $p\in
\C^+$.  It follows that \baa |K(p)|\le e^{-\a M|p+\b|^q}<1,\quad
p\in\C^+. \label{Kest}\eaa Hence $K\in H^r$ for all $r\in
[1,+\infty]$.
\begin{proposition}\label{prop1} Let $\b>0$ and a rational number $q\in(\frac{1}{2},1)$ be given.
 Let $v\in L_2(\R^+)$, $V=\L v\in H^2$. For $\a>0$, set $V_\a\defi K_{\a,\b,q} V$,
$v_\a\defi \F^{-1}V_\a(i\o)|_{\o\in\R}$. Then $V_\a\in H^2$ and
$v_\a\to v$ in $L_2(\R^+)$ as $\a\to 0$, $\a>0$.
\end{proposition}
\par
{\it Proof.} Clearly, $ V_\a(i\o)\to V(i\o)$ as $\a\to 0$ for a.e.
$\o\in\R$. By (\ref{K}), $V_\a\in H^2$. In addition,
$|K_{\a,\b,q}(i\o)|\le 1$. Hence $|V_\a(i\o)-V(i\o)|\le
2|V(i\o)|$. We have that
$\|V(i\o)\|_{L_2(\R)}=\|v\|_{L_2(\R^+)}<+\infty$.
 By Lebesgue
Dominance Theorem, it follows that \baaa \left\|V_\a(i\o)-V(i\o)
\right\|_{L_2(\R)}\to 0\quad\hbox{as}\quad \a\to 0.
 \eaaa
Hence $v_\a\to v$ in $L_2(\R^+)$ as $\a\to 0$. Then the proof
follows. $\Box$
\par
The inverse Fourier transform
$k(t)=\F^{-1}K_{\a,\b,q}(i\o)|_{\o\in\R}$ can be viewed as a
smoothing kernel; $k(t)=0$ for $t<0$. It can be seen that $k$ has
derivatives of any order.
\par
Denote by $\D$ the set of  functions $v:\R^+\to \R$ such that
there exist $\a>0$, $\b>0$, and a rational number $q\in
(\frac{1}{2},1)$,  such that $\w V\in H^{2},$ where $\w
V(p)=K_{\a,\b,q}(p)^{-1}V(p)$, $V=\L v$.
\par
The set $\D$ includes outputs of the convolution integral
operators
with the kernels $k(t)$. By Proposition \ref{prop1}, it follows that
the set $\D$ is everywhere dense in $L_2(\R^+)$.
\section{The main
result} Set $F(x,\cdot)\defi\L f(x,\cdot)$, where $x>0$ is given,
and $G_k\defi\L g_k$,  $k=0,1$.
\begin{theorem}\label{ThM}
Let the functions $f$ and $g_k$ are such that there exists $y>0$,
$\a>0$, $\b>0$, a rational number $q\in (\frac{1}{2},1)$, such
that $\w G_k\in H^{2}$, $\w F(x,\cdot)\in H^2$ for a.e. $x>0$ and
 $\int_0^y\|\w F(s,\cdot)\|_{H^2}ds<+\infty$, where \baa \w F(x,p)\defi \frac{F(x,p)}{K(p)},\qquad \w
G_k(p)\defi \frac{
 G_k(p)}{K(p)}, \label{est}\eaa and where the function $K=K_{\a,\b,q}$ is defined by
(\ref{K}) (in particular, this means that $g_k\in \D$ and
$f(x,\cdot)\in \D$ for a.e. $x\in[0,y]$). Then there exists an
unique solution $u(x,t)$ of problem
 (\ref{parab}) in the domain $(0,y)\times\R^+$ in the class
 $\W(y)$. Moreover, there exists a constant $C(y)=C(a,b,c,\a,\b,q,y)$ such
 that
\baaa   \|u\|_{\W(y)}\le C(y)\Bigl(\|\w G_1 \|_{H^2}+\|\w
G_2\|_{H^2}+\int_0^x\|\w F(s,\cdot)\|_{H^{2}}ds\Bigr). \eaaa
\end{theorem}
\begin{remark} Theorem \ref{ThM}  requires that functions $f$ and $g_k$
are smooth in $t$; in particular, they belong to $C^{\infty}$ in
$t$. However, it is not required that $f(x,t)$ is smooth in $x$.
\end{remark}
\par
{\it Proof of Theorem \ref{ThM}}. Instead of (\ref{parab}),
consider the following problems for $p\in\C^+$: \baa
&&apU(x,p)=\frac{\p^2 U}{\p x^2}(x,p)+b\frac{\p U}{\p x}(x,p)+c
U(x,p)+F(x,p),\quad x>0, \nonumber\\ &&U(0,p)\equiv G_0(p),\quad
\frac{\p U}{\p x}(0,p)\equiv G_1(p). \label{parabU}\eaa Let
$\lambda_k=\lambda_k(p)$ be the roots of the equation
$\lambda^2+b\lambda +(c-ap)=0.$  Clearly, $\lambda_{1,2}\defi
-b/2\pm \sqrt{ap+\mu}$. Recall that $\mu>0$. It follows that the
functions $(\lambda_1(p)-\lambda_2(p))^{-1}$ and
$\lambda_k(p)(\lambda_1(p)-\lambda_2(p))^{-1}$, $k=1,2$, belong to
$ H^\infty$.
\par
 For $x\in(0,y]$, the solution of (\ref{parabU}) is \baa
U(x,p)&=&
\frac{1}{\lambda_1-\lambda_2}\biggl((G_1(p)-\lambda_2G_0(p))e^{\lambda_1
x}-(G_1(p)-\lambda_1G_0(p))e^{\lambda_2x} \nonumber\\
 &-&
\int_0^x e^{\lambda_1(x-s)}F(s,p)ds+\int_0^x
e^{\lambda_2(x-s)}F(s,p)ds
 \biggr).
 \eaa
 This can be derived, for instance,  using
 Laplace transform method applied to linear ordinary differential equation
 (\ref{parabU}), and having in mind that
   \baaa \frac{1}{\lambda^2+b\lambda +c-ap}= \frac{1}{(\lambda-\lambda_1)(\lambda-\lambda_2)}=
\frac{1}{\lambda_1-\lambda_2}\left(\frac{1}{\lambda-\lambda_1}
-\frac{1}{\lambda-\lambda_2}\right),\\
\frac{\lambda}{\lambda^2+b\lambda
+c-ap}=\frac{\lambda}{(\lambda-\lambda_1)(\lambda-\lambda_2)}=
\frac{1}{\lambda_1-\lambda_2}\left(\frac{\lambda_1}{\lambda-\lambda_1}
-\frac{\lambda_2}{\lambda-\lambda_2}\right). \eaaa
\par
Let $x\in(0,y)$, $s\in[0,x]$.  The functions
$e^{(x-s)\lambda_k(p)}$, $k=1,2$, are holomorphic in $\C^+$.
\par
We have  $$
\ln|e^{(x-s)\lambda_k(p)}|=\Re((x-s)\lambda_k(p))=(x-s)\left(
-\frac{b}{2}\pm |ap+\mu|^{1/2}\cos\frac{\Arg(ap+\mu)}{2}\right),
$$ where $k=1,2,$ $p\in\C^+$. It follows that $$
|K(p)e^{(x-s)\lambda_k(p)}|\le e^{(x-s)[-b/2+|ap+\mu|^{1/2}]-\a
M|p+\b|^q},$$ $k=1,2$, $p\in\C^+$.  Similarly, $$
|K(p)e^{\lambda_k x}|\le e^{x[-b/2+|ap+\mu|^{1/2}]-\a M|p+\b|^q}.
$$
 Since $q>1/2$, it follows that $K(p)e^{\lambda_k x}\in H^r$,
 $K(p)e^{(x-s)\lambda_k(p)}\in H^r$, $pK(p)e^{\lambda_k x}\in H^r$, and
 $pK(p)\Psi_k(p)\in H^r$,
for  $r=2$ and $r=+\infty$. Moreover, we have \baaa
 &&\sup_{s\in[0,x]}\|p^me^{\lambda_k(p)
s}G_k(p)\|_{H^2}\le C_1(x)\|\ww G_k\|_{H^2},
\\
&&\sup_{s\in[0,x]}\|p^me^{\lambda_k(p) s}K(p)\|_{H^\infty}\le
C_2(x),
 \eaaa where $m=0,1$.
Hence \baaa
\sup_{x\in[0,y]}\left\|p^m\int_0^xe^{(x-s)\lambda_k}F(s,p)ds\right\|_{H^2}
\le
\sup_{x\in[0,y]}\int_0^x\left\|e^{(x-s)\lambda_k}p^mF(s,p)\right\|_{H^2}ds
\\
\le
\sup_{x\in[0,y]}\int_0^x\|p^me^{\lambda_k(x-s)}K(s)\|_{H^\infty}\|\ww
F(s,p)\|_{H^2}ds \le C_2(y)\int_0^y\|\w F(s,p)\|_{H^2}ds, \eaaa
where $m=0,1$. Here  $C_1(x)$, $C_2(x)$ are constants that depend
on $a,b,c,\a,\b,q,x$.  It follows that $p^me^{\lambda_k
x}G_m(p)\in H^2$ and
 $p^m\int_0^xe^{(x-s)\lambda_k}F(p,s)ds\in H^2$
for any $x>0$, $m=0,1$,  $k=1,2$.
\par
Recall that $\lambda_k=\lambda_k(p)$.
Let $$N\defi
\left\|\frac{1}{\lambda_1-\lambda_2}\right\|_{H^\infty}+\sum_{k=1,2}
\left\|\frac{\lambda_k}{\lambda_1-\lambda_2}\right\|_{H^\infty}.
 $$
\par
It  follows from the above estimates that \baa
\|p^mU(x,p)\|_{H^2}\le N\left(C_1(y)\sum_{k=1,2} \left\|\w
G_k\right\|_{H^2}+ C_2(y)\int_0^x\|\w
F(s,p)\|_{H^2}ds\right),\quad m=0,1. \label{est1} \eaa
\par
It follows that the corresponding inverse Fourier transforms
$u(x,\cdot)=\F^{-1}U(x,i\o)|_{\o\in \R}$, $\frac{\p u}{\p t}
(x,\cdot)=\F^{-1}(pU(x,i\o)|_{\o\in \R})$ are well defined and are
vanishing for $t<0$.
 In addition, we have that
$\overline{U(x,i\o)}=U(x,-i\o)$ (for instance,
 $\overline{K(i\o)}=K(-i\o)$,
 $\overline{e^{(x-s)\lambda_k(i\o)}}=e^{(x-s)\lambda_k(-i\o)}$, etc). It follows that
 the inverse of Fourier transform
 $
u(x,\cdot)=\F^{-1}U(x,\cdot)
 $
 is real.
\par
Further, we have that \baa \frac{\p U}{\p x}(x,p)&=&
\frac{1}{\lambda_1-\lambda_2}\biggl((G_1(p)-\lambda_2G_0(p))\lambda_1e^{\lambda_1
x}-(G_1(p)-\lambda_1G_0(p)) \lambda_2e^{\lambda_2x} \nonumber\\
 &-&
\lambda_1\int_0^x e^{\lambda_1(x-s)}F(s,p)ds+\lambda_2\int_0^x
e^{\lambda_2(x-s)}F(s,p)ds
 \biggr).
 \label{est2}\eaa
 Since $\lambda_1(p)\lambda_2(p)=c-ap$, we obtain again
that \baa \left\|\frac{\p U}{\p x}(x,p)\right\|_{H^2}\le C_3(y)
\left(\sum_{k=1,2} \left\|\w G_k\right\|_{H^2}+ \int_0^x\|\w
F(s,p)\|_{H^2}ds\right).\label{est3}\eaa By (\ref{parabU}),
$\p^2U/\p x^2$ can be expressed as a linear combination of
$F,G_k,U,pU,\p U/\p x$. By (\ref{est1})-(\ref{est3}), \baaa
\left\|\frac{\p^2 U}{\p x^2}(x,p)\right\|_{H^2}\le C_4(y) \left(
\left\|\frac{\p U}{\p x}(x,p)\right\|_{H^2}+
\sum_{m=0,1}\left\|p^mU(x,p)\right\|_{H^2}+ \|
F(x,p)\|_{H^2}\right). \eaaa
 We have that $|K(p)|<1$ on $C^+$ and $\|F(s,p)\|_{H^2}\le\|\w F(s,p)\|_{H^2}$.
 It follows that \baa \left\|\frac{\p^2 U}{\p
x^2}(x,p)\right\|_{H^2}\le C_5(y) \left(\sum_{k=1,2} \left\|\w
G_k\right\|_{H^2}+ \int_0^x\|\w
F(s,p)\|_{H^2}ds\right).\label{est4} \eaa Here $C_k(y)$ are
constants that depend on $a,b,c,\a,\b,q,y$. By
(\ref{est1})-(\ref{est4}), estimate (\ref{est}) holds.
\par
  Therefore, $u(x,\cdot)=\F^{-1}U(x,i\o)|_{\o\in \R}$ is the solution of (\ref{parab}) in $\W(y)$.
  The uniqueness is ensured by the linearity of the problem, by
  estimate (\ref{est}), and by the fact that $\L u(x,\cdot)$, $\L (\p^k u(x,\cdot)/\p x^k)$,
  and  $\L (\p u(x,\cdot/\p t)$ are well defined on $\C^+$
   for any $u\in\W(y)$.
  This completes the proof of Theorem \ref{ThM}. $\Box$
\begin{remark}
It can be seen from the proof that it is crucial that
$u(x,0)\equiv 0$. Non-zero initial conditions can  not be
included.
\end{remark}
\subsection*{References}$\hphantom{xx}$
 Duren, P. {\it Theory of $H^p$-Spaces.}
1970. Academic Press, New York.
\par
Tikhonov, A. N. and Arsenin, V. Y. {\it Solutions of Ill-posed
Problems.} 1977. W. H. Winston, Washington, D. C.
\end{document}